\documentclass[11pt]{article}
\usepackage{graphicx, amsmath, amsthm, amsfonts, enumerate, authblk,
  amssymb, graphics, fullpage, esint}
\usepackage{empheq}
\usepackage{bigints, bbm}
\usepackage{epsfig}
\usepackage{amsmath}
\usepackage{scrlfile}
\usepackage{srcltx}
\usepackage{rotating}
\usepackage{color}
\usepackage[mathscr]{euscript}
\usepackage{relsize}
\usepackage{mathrsfs} 

\newtheorem{theorem}{Theorem}[section]

\newtheorem{proposition}[theorem]{Proposition}
\newtheorem{definition}{Definition}[section]

\newtheorem{remark}{Remark}

\makeatletter
\def\blfootnote{\xdef\@thefnmark{}\@footnotetext}
\makeatother
\allowdisplaybreaks

\newcommand{\cH}{{\cal H}}

\DeclareMathOperator{\dom}{dom}
\DeclareMathOperator{\clos}{clos}
\DeclareMathOperator{\ran}{ran}

\DeclareMathOperator{\im}{Im}

\DeclareMathOperator*{\slim}{s-lim}

\usepackage{soul}
\usepackage{xcolor}
\usepackage[normalem]{ulem}

\def\e{\varepsilon}

\begin{document}
\title{Mathematical Heritage of Sergey Naboko: Functional Models of Non-Self-Adjoint Operators}
\author{Alexander V. Kiselev and Vladimir Ryzhov}

\maketitle

\section{Dilation theory for dissipative operators}
\label{sec:functional-model}


Functional model construction for a contractive linear operator~$T$
acting on a Hilbert space~$K$ is a well developed domain of the
operator theory.
Since  pioneering works by B.~Sz.-Nagy, C.~Foia\c{s}~\cite{MR2760647},
P.~D.~Lax, R.~S.~Phillips~\cite{MR0217440},  L.~de~Branges,
J.~Rovnyak~\cite{deBranges, deBrangesRovnyak}, and
M.~Liv\v{s}ic~\cite{Livsic1954}, this research area attracted many
specialists in operator theory, complex analysis, system control,
gaussian processes and other disciplines.
Multiple studies culminated in the development of a comprehensive theory
complemented by various applications,
see~\cite{DyM,Fu,Nikolski,NikolskiiKhrushchev, NikolskiiVasyunin1998}
and references therein.
%

The underlying idea of functional model is the fundamental theorem
of B.~Sz.-Nagy and C.~Foias stating that for a dissipative operator $L$
under the assumption~$ \mathbb C_- \subset \rho(L)$
(dissipative operators satisfying this condition are
called {\em maximal}),
there exists a selfadjoint dilation of~$L$, which  is a selfadjoint
operator~$\mathscr L$ on a wider
space~$\mathcal H \supset K$
such that
\begin{equation}\label{DilationCondition}
(L-zI)^{-1} = P_K (\mathscr L - zI)^{-1}|_K, \quad
z \in \mathbb C_-,
\end{equation}
where~$P_K$ is an orthogonal projection
from $\mathcal H$ onto $K$.

{In applications}, such a
dilation~$\mathscr L$  should
be {\it minimal};
it should not contain any
reducing selfadjoint parts not
related to the operator~$L$.
Mathematically the minimality condition is
expressed as the equality
\[
\clos\bigvee_{z \notin R} {(\mathscr L -  zI)^{-1}
\mid_K }  =
\mathcal H,
\]
where $\mathcal H$ is the dilation space~$\mathcal H
\supset  K$.
Construction of a dilation satisfying this
 condition is a non-trivial task successfully solved
for contractions
by Sz.-Nagy and Foia\c{s}~\cite{MR2760647}
with the help of Neumark's theorem~\cite{Naimark1943},
and by B.~Pavov~\cite{MR0510053, Drogobych}
for two
important cases of dissipative
operators arising in
mathematical physics and
 successfully
extended later to a general setting
(more on this in the following sections).
%


The functional model theory of
non-selfadjoint operators
studies operators~$L$ which have
no non-trivial reducing
selfadjoint parts.
Such operators are called {\em completely
non-selfadjoint} or, using
a less accurate term, {\em simple}.
{In what follows,
all  non-selfadjoint operators are
assumed closed, densely defined and simple,
with regular points in both lower
and upper half-planes.}

\subsection{Additive perturbations}

Let $A = A^*$ be a selfadjoint unbounded operator on a Hilbert space~$K$ and
$V$ a bounded (for simplicity) non-negative operator~$V  = V^* = \alpha^2/2\geq 0$, where
$\alpha = (2V)^{1/2}$. Let $L = A + \frac{i}{2} \alpha^2$. The operators~$A$ and~$V = \alpha^2/2$ are the
real and imaginary parts
of~$L$ defined on $\dom(L) = \dom(A)$.

Following Pavlov, denote $E = \clos \ran \alpha$ and define the dilation
space as the direct sum of
$K$ and the equivalents of incoming and
outgoing channels
{of the Lax-Phillips scattering theory, see \cite{MR0217440},}~ $\mathcal D_\pm = L^2(\mathbb
R_\pm, E)$,
\begin{equation}\label{DilationSpace}
\mathcal H =  \mathcal D_- \oplus K \oplus \mathcal D_+.
\end{equation}
Elements of~$\mathcal H$ are represented as
three-component vectors~$(v_-, u,
v_+)$ with
$v_\pm \in \mathcal D_\pm$ and $u \in K$.
The action of~$\mathscr L$ on the channels~$\mathcal D_\pm$
is defined by
$\mathscr L : (v_-,0, v_+) \mapsto ( iv_-^\prime,0,iv_+^\prime)$.
The self-adjointness of~$\mathscr L = \mathscr L^*$ and
the requirement~(\ref{DilationCondition}) lead to
the form of dilation~$\mathscr L$
suggested in~\cite{MR0510053},
\begin{equation}\label{DilationForm}
 \mathscr L  \left(
 \begin{array}{l} v_-
                     \\ u
                     \\   v_+
\end{array}\!\!\right)
=
 \begin{pmatrix} i\frac{dv_-}{dx}
                     \\ A u + \frac{\alpha}{2} \left[v_+(0) + v_-(0)  \right]
                     \\   i\frac{dv_+}{dx}
\end{pmatrix},
\end{equation}
defined on the domain
\[
 \dom(\mathscr L) = \left \{(v_-, u, v_+ ) \in \mathcal H \mid
 v_\pm \in W_2^1(\mathbb R_\pm, E), u \in \dom(A),
 v_+(0) - v_-(0) = i\alpha u
 \right\}
\]
The ``boundary condition''~$v_+(0) - v_-(0) = i\alpha u$
can be interpreted
as a coupling between the incoming and outgoing channels
$\mathcal D_\pm$, realised by the
imaginary part of~$L$ acting
on $E$.
The characteristic function of $L$ is the contractive
operator-valued
function defined by the formula
\begin{equation}\label{DissCharF}
S(z) =  I_E + i\alpha(L^* - zI)^{-1}\alpha : E \to E,
\quad z \in \mathbb C_+.
\end{equation}
%


Owing to the general theory~\cite{MR2760647},
the operator~$L$ is unitary equivalent
to its model 
in the spectral representation of~$\mathscr L$
in accordance with~(\ref{DilationCondition}).

Due to the operator version of Fatou's
theorem~\cite{MR2760647},
non-tangential boundary values of the
function~$S$
exist in the strong operator topology
almost everywhere on the real line.
Denote~$S = S(k) =  \slim_{\varepsilon\downarrow 0}
S(k + i\varepsilon)$, a.~e. $k \in \mathbb R$.
Similarly, let~$S^* = S^*(k)
:= \slim_{\varepsilon\downarrow 0}[S(k + i\varepsilon)]^*$,
which exists for almost all~$k \in \mathbb R$.
%
%
The symmetric form of the functional model
is obtained
by factorisation and completion
of the dense linear set of
vector-valued functions from the
space~$L^2(E) \oplus L^2(E)$
with respect to the norm
\begin{equation}\label{ModelSpace}
 \left\|\binom{\tilde g}{g}\right\|_{\mathscr H}^2 :=
 \int\limits_{\mathbb R}
  \left\langle
    \left(
      \begin{array}{cc} I &  S^* \\  S & I \end{array}
    \right)
    \binom{\tilde g}{g},
    \binom{\tilde g}{g}
  \right\rangle_{ E\oplus  E}
 dk
\end{equation}
%
Note that the
elements of~$\mathscr H$ are not individual functions from
$L^2(E) \oplus L^2(E)$, but rather equivalence classes
formed after factorization over elements with zero
$\mathscr H$\nobreakdash--norm, {followed by}
completion~\cite{NikolskiiKhrushchev,
NikolskiiVasyunin1989}.
It is easily seen that for each $\binom{\tilde g}{g} \in \mathscr H$
the expressions
$g_+:= S\tilde g + g$ and $g_-:=\tilde g +  S^* g$ are
in fact usual square summable
vector-functions
from~$L_2( E)$.

The space~$\mathscr H = L_2\begin{psmallmatrix}
                            I & S^* \\ S & I
                           \end{psmallmatrix}
$ with the norm defined
by~(\ref{ModelSpace}) turns out to be the spectral representation space of the self-adjoint
dilation~$\mathscr L$ of the operator~$L$. Henceforth we will denote the corresponding unitary mapping of $\mathcal{H}$ onto $\mathscr{H}$ by $\Phi$.
It means that
the operator of multiplication by the independent variable
acting on~$\mathscr H$, i.e., the
operator~$f(k) \mapsto k f(k), $ is unitary equivalent
to the dilation~$\mathscr L$.
Hence, for~$z \in\mathbb
C \setminus \mathbb R$, the
mapping~$\binom{\tilde g}{g} \mapsto (k - z)^{-1}
\binom{\tilde g}{g}$, where $\binom{\tilde g}{g} \in \mathscr
H$
is unitary equivalent to the resolvent of~$\mathscr L$
and therefore~$L$ is mapped to its
functional model (with the symbol $\simeq$ denoting unitary equivalence),
\begin{equation}\label{ModelInModelSpace}
 (L - zI)^{-1} \simeq\left . P_{\mathscr K} (k - z)^{-1} \right
 |_{\mathscr K}, \quad z \in \mathbb C_-
\end{equation}

The incoming and outgoing subspaces of the dilation
space~$\mathscr H$ admit the form
\begin{equation*}
 {\mathscr D}_+ := \binom{H_2^+(E)}{0}, \quad
 {\mathscr D}_- := \binom{0}{H_2^-( E)}, \quad
 \mathscr K := \mathscr H \ominus
  \left[
    {\mathscr D}_+ \oplus {\mathscr D}_-
  \right]
\end{equation*}
where $H_2^\pm(E)$ are the Hardy classes of $
E$-valued vector-functions analytic in $\mathbb C_\pm$
{and $\mathscr D_\pm = \Phi\mathcal D_\pm$}.
As usual~\cite{MR822228}, the functions from vector-valued Hardy
classes~$H_2^\pm(E)$ are identified with their
boundary values existing almost everywhere
on the real line.
They form two complementary mutually orthogonal subspaces,
so that
$L_2(E) = H_2^+(E)\oplus H_2^-(E)$.

The image~$\mathscr K$ of $K$ under the spectral
mapping $\Phi$ of the dilation space~$\mathcal H$
to~$\mathscr H$ is the subspace
\begin{equation*}
 \mathscr K =
   \left\{
     \binom{\tilde g}{g} \in \mathscr H \,\, \mid\,\,
       \tilde g +  S^* g \in H_2^-( E), \,
        S \tilde g + g \in H_2^+ ( E)
   \right\}
\end{equation*}
The orthogonal projection~$P_\mathscr K$ from
$\mathscr H$ onto
$\mathscr K$ is defined by formula~(\ref{2.1.3.2}).
Note that the following definition
has to be understood
on the dense set of functions from
$L^2(E) \oplus L^2(E)$ in $\mathscr H$.
\begin{equation}\label{2.1.3.2}
  P_{\mathscr K} \binom{\tilde g}{g} =
    \binom{\tilde g - P_+ ( \tilde g +   S^* g)}
          {       g - P_- ( S \tilde g +   g)},
    \quad \tilde{g} \in L_2( E), \; g \in L_2( E),
\end{equation}
where $P_\pm$ are the orthogonal projections from $L_2$
onto the Hardy
classes~$H_2^\pm$.


%

%
%
%
%

\section{Naboko's  functional model for a family of additive perturbations}
\label{sec:nabok-spectr-form}

The model approach to the analysis of dissipative operators 
outlined above relies exclusively on the knowledge of a characteristic
function of a dissipative completely non-selfadjoint operator~$L$.
The properties of the operator are
expressed in terms of its characteristic function, i.~e.,
in the language of analytic operator-valued functions theory.
This represents the true value of the functional model approach:
all the abstract results obtained using model techniques
become immediately available, once the
characteristic function of the operator is known.
%
%
%
%
%
%

Successful applications of the functional
model approach for
contractions and dissipative operators {have} inspired the
search for models of
non-dissipative operators.
The attempts to follow the blueprints of Sz.-Nagy-Foias
and Lax-Philitps meet serious challenges
rooted in the absence of a self-adjoint dilation for
such operators.
\smallskip

The breakthrough came in the late seventies with the
publication of papers~\cite{MR0500225,nabokozapiski2} and
especially~\cite{MR573902} by S.~Naboko, who
found a way to represent a non-dissipative
operator in a model space of a suitably
chosen dissipative one.
Apart from the model construction, his works
largely contributed to the development of various areas
in the non-self-adjoint operator theory.
In contrast to the earlier results, {his} model
representation  does not rely on the
uniqueness (up to a unitary equivalence) of the characteristic function
of a completely non-selfadjont operator.
Based on the dilation~(\ref{DilationForm}), the
paper~\cite{MR573902} provides an isometry
between the dilation space~(\ref{DilationSpace}) and
the model space~(\ref{ModelSpace}) in an
explicit form.
This explicitness plays a crucial r\^ole in passage to the
model representation for non-dissipative operators
using nothing more than Hilbert
resolvent identities.
All the building blocks of the method are clearly presented in
terms of the original problem, which is especially appealing from the
applications' perspective.
We next give a brief overview of the key ideas presented in~\cite{MR0500225,nabokozapiski2,MR573902}.

\subsection{Isometric map between the dilation and model spaces}\label{sec:nsa_setup}

Consider a non-self-adjoint operator
\begin{equation}\label{eqn:AV}
 L = A + iV
\end{equation}
acting in the Hilbert space~$K$, where~$A= A^*$
 and~$V= V^*$
is~$A$-bounded with the relative bound less than $1$.
The domains of~$A$ and $L$ coincide and the
operator~$L$ is closed.
Note that $V$ can be written in the
form~$V = \frac{\alpha J \alpha}{2}$ with
$\alpha = \sqrt{2|V|}$,
$J:= \mathop{sign} V : E\to E$ defined
according to the functional calculus of
self-adjoint operators.
Like in~(\ref{DissCharF}), $E := \clos\ran \alpha$.
The characteristic function of~$L$ admits
the form~(see, e.g.,~\cite{Strauss1960})
\begin{equation}\label{eqn:AVCharF}
 \Theta(z)= I_E + i J\alpha(L^* -zI)^{-1} \alpha : E \to E,
 \quad z \in \rho(L^*).
\end{equation}

Alongside with $L$
introduce the operator~${L^{\phantom{s}}\!\!}^{\scriptscriptstyle{||}}$
on the same
domain~$\dom({L^{\phantom{s}}\!\!}^{\scriptscriptstyle{||}}) =
\dom(L)$ as follows:
\begin{equation}\label{eqn:AVMod}
 {L^{\phantom{s}}\!\!}^{\scriptscriptstyle{||}} := A + i|V| = A + i
\frac{\alpha^2}{2}.
\end{equation}
The operator~${L^{\phantom{s}}\!\!}^{\scriptscriptstyle{||}}$
is precisely the dissipative operator of the preceding Section.
The work~\cite{MR573902} contains the model construction,
the definition of the isometry~$\Phi:
\mathcal H \to \mathscr H$ from (a dense set in) the
dilation space~(\ref{DilationSpace}) to the model space~(\ref{ModelSpace})
of ${L^{\phantom{s}}\!\!}^{\scriptscriptstyle{||}}$,
which is a preliminary step towards the model
for its additive perturbations
of the form~(\ref{eqn:AV}).
Note that the characteristic function~$S$
of~${L^{\phantom{s}}\!\!}^{\scriptscriptstyle{||}}$
is given by the expression~(\ref{DissCharF}) where $L$ is replaced by
${L^{\phantom{s}}\!\!}^{\scriptscriptstyle{||}}$:
\begin{equation}\label{eqn:AVmodCharF}
 S(z)= I_E + i \alpha( {L^{\phantom{s}}\!\!\!}^{\scriptscriptstyle{-||} }
-zI)^{-1} \alpha,
 \quad z \in \rho({L^{\phantom{s}}\!\!\!}^{\scriptscriptstyle{-||} }), \quad {L^{\phantom{s}}\!\!\!}^{\scriptscriptstyle{-||} }:=({L^{\phantom{s}}\!\!}^{\scriptscriptstyle{||} })^*.
\end{equation}
%

The argument of~\cite{MR0500225}  shows that the
characteristic functions of $L$ and
${L^{\phantom{s}}\!\!}^{\scriptscriptstyle{||}}$ are
related via the Potapov-Ginzburg
operator linear-fractional  transformation,
or PG-trans\-form~\cite{AzizovIokhvidov}.
This fact is essentially geometric.
It relates contractions on Kre\u{\i}n spaces
(i.~e., the spaces with an indefinite metric defined
by the involution~$J = J^* = J^{-1}$) to contractions
on Hilbert spaces.
The PG-trans\-form is invertible and
the following assertion pointed out in~\cite{MR0500225}
holds.
\begin{proposition}
The characteristic function~(\ref{eqn:AVCharF})
of~$L = A + iV$ is $J$\nobreakdash-contractive on its domain
and the PG\nobreakdash-transform
maps $\Theta$ to the contractive characteristic
function of~${L^{\phantom{s}}\!\!}^{\scriptscriptstyle{||}} = A + i|V|$
defined by~(\ref{eqn:AVmodCharF}), as follows:
\begin{equation}\label{PGtransform}
\Theta \mapsto S = - (\chi^+ - \Theta \chi^-)^{-1}(\chi^- - \Theta \chi^+),
\qquad
S \mapsto \Theta = (\chi^- + \chi^+ S)(\chi^+ + \chi^- S)^{-1},
\end{equation}
where $\chi^\pm = \frac{1}{2} (I_E \pm J)$ are orthogonal projections onto
the subspaces $\chi^+ E$ ($\chi^- E$, respectively).
\end{proposition}
\noindent

It appears somewhat unexpected that two operator-valued functions
connected by formulae~(\ref{PGtransform})   can be explicitly
written down in terms of their ``main
operators''~$L$ and $ {L^{\phantom{s}}\!\!\!}^{\scriptscriptstyle{-||} }$.
This relationship between the characteristic
functions of~$L$ and
${L^{\phantom{s}}\!\!}^{\scriptscriptstyle{||}}$  goes in fact
much deeper, see~\cite{Arov,AzizovIokhvidov}.
In particular, the self-adjoint dilation
of~${L^{\phantom{s}}\!\!}^{\scriptscriptstyle{||}}$
and the $J$--self-adjoint dilation of $L$ are also related via a suitably
adjusted version of the PG\nobreakdash-transform.
Similar statements hold for the corresponding
linear systems or ``generating operators'' of the
functions~$\Theta$ and $S$, see~\cite{Arov,AzizovIokhvidov}.
This fact is crucial for the construction of a model of a
general closed, densely defined
non-self-adjoint operator,
see~\cite{Ryzhov_closed}.

\smallskip

Assume as usual that~the
operator~${L^{\phantom{s}}\!\!}^{\scriptscriptstyle{||}}$
is completely non-self-adjoint, and
let~$\mathscr L $ be the minimal self-adjoint dilation of
${L^{\phantom{s}}\!\!}^{\scriptscriptstyle{||}}$ of the
form~(\ref{DilationForm}).
\begin{theorem}[\cite{MR573902}, Theorem 2]\label{ModelTheoremDiss}
There exists a  mapping $\Phi$ from the
dilation space~$\mathcal H$ onto  Pavlov's model space~$\mathscr H$
defined by~(\ref{ModelSpace})
with the following properties.
\begin{enumerate}
  \item  
    $\Phi$ is  isometric.
  \item  
    $\tilde g +  S^* g = \mathscr F_+  h$,
    $  S \tilde g + g = \mathscr F_-  h$, where
    $\binom{\tilde g}{g} = \Phi  h$, $ h \in \mathcal H$
  \item
    $
   \Phi \circ (\mathscr L - zI )^{-1} = (k - z)^{-1} \circ \Phi,
   \quad z \in \mathbb C \setminus\mathbb R
    $
  \item
    $ \Phi\mathcal H  = \mathscr H$, \quad $\Phi \mathcal D_\pm = {\mathscr
D}_\pm$, \quad $\Phi  K = \mathscr K$
  \item
    $\Phi \circ (\mathscr L - zI)^{-1} = (k - z)^{-1} \circ
    \Phi,
     \quad z\in \mathbb C \setminus \mathbb R$.
\end{enumerate}
Here the bounded maps~$\mathscr F_\pm : \mathcal H \to L^2(\mathbb R, E)$
are defined by the formulae
 \begin{align*}
  \mathscr F_+ : h &  \mapsto - \frac{1}{\sqrt{2\pi}}\,
 \alpha ( {L^{\phantom{s}}\!\!}^{\scriptscriptstyle{||}} - k + i0)^{-1} u +
 S^*(k) {\hat{v}}_-(k) +    \hat{v}_+(k),
  \\
  \mathscr F_- : h &  \mapsto - \frac{1}{\sqrt{2\pi}}\,
   \alpha ( {L^{\phantom{s}}\!\!\!}^{\scriptscriptstyle{-||}}  - k - i0)^{-1}
u + {\hat{v}}_-(k) + S(k) \hat{v}_+(k),
 \end{align*}
where $h = (v_-, u, v_+) \in \mathcal H$ and $\hat v_\pm$ are the Fourier
transforms of $v_\pm \in L^2(\mathbb R_\pm, E)$.
%
\end{theorem}

\subsection{Model representation of additive perturbations}\label{sec:additive}

Theorem~\ref{ModelTheoremDiss} opens a possibility of
expressing a larger class of perturbations of~$A$
in the model space~$\mathscr H$.
Namely, consider operators in~$K$  of the form
\begin{equation}\label{operator:Lkappa}
 L^\varkappa = A + \frac{\alpha \varkappa\alpha} {2},  \quad
 \dom(L^\varkappa) = \dom(A),
\end{equation}
where $\varkappa$ is a bounded operator in $E$.
The family $\{L^\varkappa \mid {\varkappa: E\to E} \}$
includes $A$ for $\varkappa = 0$, the dissipative
operator~${L^{\phantom{s}}\!\!}^{\scriptscriptstyle{||} }$
for $\varkappa = iI_E$,
its adjoint
${L^{\phantom{s}}\!\!\!}^{\scriptscriptstyle{-||} }$  for
$\varkappa = -i I_E$, as well as self-adjoint and
non-self-adjoint operators corresponding to other values of
the ``parameter''~$\varkappa$.
In particular, the non-dissipative operator~$L = A + iV = A + i\frac{\alpha
J \alpha}{2}$ of~(\ref{eqn:AV}) is recovered by
putting~$\varkappa = iJ$.
Representations of the
resolvent~$(L^\varkappa - zI)^{-1}$, $z \in \rho(L^\varkappa)$
in the model space~$\mathscr H$ are obtained using the
properties of $\mathscr F_\pm$ given in
Theorem~\ref{ModelTheoremDiss} and resolvent identities
for~$({L^{\phantom{s}}\!\!}^{\scriptscriptstyle{||} } - zI)^{-1}$,
$({L^{\phantom{s}}\!\!\!}^{\scriptscriptstyle{-||} } - zI)^{-1}$,
and
$(L^\varkappa - zI)^{-1}$.
The key component of the proofs is the representation of
$\mathscr F_\pm (L^\varkappa - zI)^{-1} u $ in terms of $\mathscr
F_\pm u$ for $u \in K$.
For instance,  it can be shown that
there exist two analytic
operator-functions~$\Theta_\varkappa^\prime, \Theta_\varkappa :E \to E$,
bounded in $\mathbb C_-$,  $\mathbb C_+$ respectively,
such that
for $z_0 \in \rho(L^\varkappa)$, $\im
z_0 < 0$, and all $u \in K$
\begin{equation}\label{eqn:FpmEquations}
\begin{aligned}
 \mathscr F_+ (L^\varkappa
- z_0I)^{-1} u
 &= \frac {1}{k - z_0 }
(\mathscr F_+ u)(k-i0)
-\frac {1}{k - z_0 }
\Theta_\varkappa^\prime(k-i0)
 [ \Theta_\varkappa^\prime(z_0)]^{-1}
 (\mathscr F_+ u)(z_0)
 \\
 \mathscr F_- (L^\varkappa
- z_0I)^{-1} u
 &= \frac {1}{k - z_0}
 (\mathscr F_- u)(k+i0)
-\frac {1}{k - z_0 }
\Theta_\varkappa(k+i0)
 [ \Theta_\varkappa^\prime(z_0)]^{-1}
 (\mathscr F_+ u)(z_0)
\end{aligned}
\end{equation}
Here  $\mathscr F_\pm u  \in H_2^\mp(E)$ since $u \in K$
and $(\mathscr F_+u)(z_0) = (\tilde g + S^*g)(z_0)$ is the
analytic continuation of the function~$(\tilde g + S^*g)$ to the
point~$z_0$ in the lower half-plane.
The possibiliy to
express $\mathscr F_\pm (L^\varkappa - z_0I)^{-1} u$
using the spectral mappings~$\mathscr F_\pm$
applied to $u\in K$ found
on the right hand side of~(\ref{eqn:FpmEquations})
is the key ingredient of calculations
leading to the main theorem.

\begin{theorem}[Model Theorem,~\cite{MR573902}]\label{ModelTheoremDiss2}
 If $z_0 \in \mathbb C_- \cap \rho (L^\varkappa)$
 and $\binom{\tilde g}{ g} \in  \mathscr K$, then
 \[
  \Phi (L^\varkappa - z_0 I)^{-1} \Phi^*\binom{\tilde g}{ g} =
  P_{\mathscr K} \frac{1}{k - z_0}
  \begin{pmatrix}
   \tilde g \\
   g - \frac{1 + i\varkappa}{2 }
\left[\Theta_\varkappa^\prime\right(z_0)]^{-1}(\tilde g +S^*g)(z_0)
  \end{pmatrix}
 \]
If $z_0 \in \mathbb C_+ \cap \rho (L^\varkappa)$
 and $\binom{\tilde g}{ g} \in K$, then
 \[
  \Phi (L^\varkappa - z_0 I)^{-1} \Phi^*\binom{\tilde g}{ g} =
  P_ {\mathscr K} \frac{1}{k - z_0}
  \begin{pmatrix}
   \tilde g -
    \frac{1 - i\varkappa}{2 }
\left[\Theta_\varkappa\right(z_0)]^{-1}(S\tilde g +g)(z_0)
   \\
   g
  \end{pmatrix}
 \]
\end{theorem}
\noindent

\subsection{Smooth vectors and the absolutely continuous subspace}\label{sec:smooth}

In  \cite{MR573902, MR1252228} Sergey Naboko introduced absolutely continuous subspaces of the family $L^\varkappa$. He always admired Mark Kre\u\i n, and in particular liked to quote him as saying: ``the major instruments of
self-adjoint spectral analysis arise from the Hilbert space geometry,
whereas in the non-self-adjoint setup the modern complex analysis has to take
the role of the main tool''. It is therefore not surprising that his definition of spectral subspaces is formulated in the language of complex analysis.

In the functional model space $\mathscr{H}$
consider two subspaces $\mathscr{N}^\varkappa_\pm$ defined as follows:
\begin{equation*}
   \mathscr{N}^\varkappa_\pm:=\left\{\binom{\widetilde{g}}{g}\in\mathscr{H}:
  P_\pm\left(\chi_\varkappa^+(\widetilde{g}+S^*g)+\chi_\varkappa^-(S\widetilde{g}
  +g)\right)=0\right\}, \text{ where }
  \chi_\varkappa^\pm:=\frac{I\pm{\rm i}\varkappa}{2}.
\end{equation*}

These subspaces are then characterised in terms of the resolvent of the operator $L^\varkappa.$ This, again, can be seen as a consequence of a much more general argument (see, {\it e.g.}, \cite{Ryzhov_closed, Ryzh_ac_sing}). 
Consider the counterparts of $\mathscr{N}^\varkappa_\pm$ in the original
Hilbert space $K:$
\begin{equation*}
  \widetilde{N}_\pm^\varkappa:=\Phi^*P_{\mathscr K}\mathscr{N}^\varkappa_\pm\,,\quad
  N_\pm^\varkappa:= \clos  \widetilde{N}_\pm^\varkappa.
\end{equation*}
Now introduce the set
$
  \widetilde{N}_{\rm e}^\varkappa:=\widetilde{N}_+^\varkappa\cap
 \widetilde{N}_-^\varkappa
$
of so-called \emph{smooth vectors} and its closure $N_{\rm e}^\varkappa (L^\varkappa):=\clos(\widetilde{N}_{\rm e}^\varkappa).$



The next assertion has been always singled out by S. Naboko in his lectures on functional models as ``the main result of the whole lecture course''. In particular, it motivates the term ``the set of smooth vectors'' used for $\widetilde{N}^\varkappa_{\rm e}$ and opens up a possibility to construct a rich functional calculus of the absolutely continuous ``part'' of the operator, leading in particular to the scattering theory (see details in the next Section).
\begin{theorem}
  \label{lem:on-smooth-vectors-other-form}
The sets $\widetilde{N}_\pm^\varkappa$ are described as follows:
\begin{equation*}
\widetilde{N}_\pm^\varkappa=\{u\in\cH: \alpha(L^{\varkappa}-z
I)^{-1}u\in H^2_\pm(E)\}.
\end{equation*}
Moreover, for the functional model image of $\tilde{N}^\varkappa_{\rm e}$ the following representation holds:
\begin{align}
&\Phi\widetilde{N}^\varkappa_{\rm e}=\biggl\{P_{\mathscr K}\binom{\widetilde{g}}{g}\in\mathscr{H}:\nonumber
\\
&\binom{\widetilde{g}}{g}\in\mathscr{H}\ {\rm satisfies}\
  \Phi(L^{\varkappa}-z I)^{-1}\Phi^*P_{\mathscr K}\binom{\widetilde{g}}{g}=
P_{\mathscr K}\frac{1}{\cdot-z}\binom{\widetilde{g}}{g}\ \ \ \forall\,z\in{\mathbb C}_-\cup{\mathbb C}_+\biggr\}.
\label{New_Representation}
\end{align}

\end{theorem}



The above Theorem together with Theorem \ref{thm:on-smooth-vectors-a.c.equality} motivated generalising the notion of the absolutely continuous subspace $\cH_{\rm ac}(L^{\varkappa})$ to the case of non-self-adjoint operators $L^\varkappa$ by identifying it with the set $N^\varkappa_{\rm e}.$
\begin{definition}
\label{abs_cont_subspace}
For a non-self-adjoint $L^\varkappa$ the absolutely continuous subspace $\cH_{\rm ac}(L^{\varkappa})$  is defined by the formula $\cH_{\rm ac}(L^{\varkappa})=N^\varkappa_{\rm e}(L^\varkappa).$

In the case of a self-adjoint operator $L^\varkappa$, $\cH_{\rm ac}(L^{\varkappa})$ is to be understood in the sense of the classical definition of the absolutely continuous subspace of a self-adjoint operator.
\end{definition}

\begin{theorem}
  \label{thm:on-smooth-vectors-a.c.equality}
  Assume that $\varkappa=\varkappa^*$
  and let  $\alpha(L^{\varkappa}-z I)^{-1}$ be a Hilbert-Schmidt
  operator for at least one point $z\in\rho(L^\varkappa)$. Then the definition $\cH_{\rm ac}(L^{\varkappa})=N^\varkappa_{\rm e}$ is
  equivalent to the classical definition of the absolutely continuous subspace of a self-adjoint operator, {\it i.e.},
  $
    N_{\rm e}^\varkappa = \cH_{\rm ac}(L^{\varkappa}).
  $
\end{theorem}

\begin{remark}
Alternative conditions, which are even less restrictive in general, that guarantee the validity of the assertion of Theorem
\ref{thm:on-smooth-vectors-a.c.equality} were obtained in \cite{MR1252228}. The
absolutely continuous subspace of a non-self-adjoint operator also admits
different definitions{~\cite{Ryzh_ac_sing}}, which in generic case can be
not equivalent to the one given above. This question is treated in full details
by Romanov in \cite{Romanov}.
\end{remark}

\subsection{Scattering theory}

The intrinsic relationship between the scattering theory and the theory of
dilations and functional models is due to \cite{MR0217440}. The fact that the
characteristic function of an arbitrary dissipative operator $L$ can be realised
as the scattering matrix of its dilation $\mathscr L$ was observed by Adamyan
and Arov in \cite{MR0206711}. This fact, as was reiterated by Sergey on many
occasions, together with Birman's seminal works on the mathematical scattering
theory, motivated his work on the construction of wave and scattering operators
in the functional model representation. With the introduction of smooth vector
sets which are dense in absolutely continuous subspaces of operators $L^\varkappa$,
it was natural to define (see \cite{MR0500225, MR573902}) the action of
exponential groups $\exp{iL^\varkappa t}$ in $\mathscr H$ as multiplication by
{$\exp(i k t)$} on the smooth vectors.

In view of the classical definition of the wave operator of a pair of self-adjoint operators,
\begin{equation*}
  W_\pm(L^0,L^\varkappa):=\slim_{t\to\pm\infty}e^{iL^0t}e^{-iL^\varkappa t}P_{\rm ac}^\varkappa,
\end{equation*}
where $P_{\rm ac}^\varkappa$ is the projection onto the absolutely continuous subspace of $L^\varkappa,$ he observed that, at least formally, for $\Phi^*P_{\mathscr K}\binom{\widetilde{g}}{g}\in\widetilde{N}_{\rm e}^\varkappa$ one has
\begin{equation}
  \label{eq:formula-0-kappaw-}
  W_-(L^0,L^\varkappa)\Phi^*P_{\mathscr K}\binom{\widetilde{g}}{g}=\Phi^*P_{\mathscr K}\binom{-(I+S)^{-1}(I+S^*)g}{g}\,,
\end{equation}
and similar formulae hold for $W_+(L^0,L^\varkappa)$, $W_\pm(L^\varkappa,L^0)$.

The need to attribute rigorous meaning to the right hand side of the latter equality, and thus to prove the existence and completeness of wave operators, motivated Sergey to investigate the boundary behaviour of operator-valued $R-$functions, see \cite{MR1036844,MR1252228} and references therein. This research has since found numerous applications in as seemingly unrelated areas as, say, the theory of Anderson localisation of stochastic differential operators. In the scattering theory (see \cite{MR1252228}) it has allowed him to prove the classical Kre\u\i n--Birman--Kuroda theorem, the invariance principle and their non-self-adjoint generalisations by following the approach sketched above. It is worth mentioning that the latter effectively blends together non-stationary, stationary and smooth formulations of the self-adjoint scattering theory.

\subsection{Singular spectrum of non-self-adjoint operators}

A major thrust of Sergey's research was towards the analysis of singular
spectral subspaces of non-self-adjoint operators. In the present section, we
mention some of his results obtained in this direction. The notation throughout
is as in Sections \ref{sec:smooth} and \ref{sec:additive}, with $\varkappa$ set
to be equal to $iJ$ with an involution $J$ (see Section \ref{sec:nsa_setup}). To
simplify the notation, we therefore consistently drop the corresponding
superscripts, as in $L=L^\varkappa$. It is further assumed throughout that the
non-real spectrum of $L$ is countable, with finite multiplicity. This latter
condition holds in particular when the perturbation $V$ is in trace
class, which we will assume satisfied ({similar results under} less
restrictive conditions are also available).

The singular subspace of $L$ is defined as follows: $N_i(L):=H\ominus N_e(L^*)$. For the operator $L^*$, it is set by $N_i(L^*):=H\ominus N_e(L)$. These definitions prove to be consistent with the classical one for self-adjoint operators due to the characterisation
$$
N_i(L)=\{u\in K: \langle ((L-t-i\e)^{-1}-(L-t+i\e)^{-1})u,v \rangle\ \to 0 \text{ as } \e\to 0 \text{ for all } v\in K\}.
$$
Define
\begin{equation}\label{eq:thetas}
\Theta_1(z)=\chi^-+S(z)\chi^+, \quad \Theta_2(z)=\chi^++S(z)\chi^-,
\quad
\Theta'_1(z)=\chi^-+S^*(\bar z)\chi^+, \quad \Theta'_2(z)=\chi^++S^*(\bar z)\chi^-,
\end{equation}
so that for the characteristic function $\Theta(z)$ one has (cf. \eqref{PGtransform})
$$
\Theta(z)=\Theta'^{*}_1(\bar z)(\Theta'^*_2)^{-1}(\bar z), \quad z\in\mathbb{C}_+; \quad
\Theta(z)=\Theta_2^*(\bar z) (\Theta_1^*)^{-1}(\bar z), \quad z\in \mathbb{C}_-.
$$
Set
$$
\widetilde N_+^i(L)=\Phi^*P_{\mathscr K} \binom{H_2^-(E)\ominus \Theta'_1 H_2^-(E))}{0},\quad
\widetilde N_-^i(L)=\Phi^*P_{\mathscr K} \binom{0}{H_2^+(E)\ominus \Theta_2 H_2^+(E))}
$$
for the operator $L$ and similarly
$$
\widetilde N_+^i(L^*)=\Phi^*P_{\mathscr K} \binom{H_2^-(E)\ominus \Theta'_2 H_2^-(E))}{0},\quad
\widetilde N_-^i(L^*)=\Phi^*P_{\mathscr K} \binom{0}{H_2^+(E)\ominus \Theta_1 H_2^+(E))}
$$
for the operator $L^*$.
The respective closures of these sets $N_+^i(L),N_-^i(L),N_+^i(L^*)$ and $N_-^i(L^*)$ are introduced in \cite{Naboko_singular}. These subspaces are invariant with respect to the resolvents of $(L-z)^{-1}$, $(L^*-z)^{-1}$. It is shown that $N_\pm^i (L)$ can be seen as spectral for $L$, representing the parts of the singular spectrum pertaining to the (closed) upper and lower half-planes, respectively. In particular, eigenvectors and root vectors of the operator $L$, corresponding to $z\in\mathbb{C}_+$ ($z\in\mathbb{C}_-$), belong to $\widetilde N_+^i(L)$ ($\widetilde N_-^i(L)$, respectively). The paper \cite{Naboko_singular} discusses the conditions of separability of
spectral subspaces under the additional condition
\begin{equation}\label{eq:separability}
\sup_{\im z>0} \max\{\|\chi^+ S(z)\chi^-\|,\|\chi^- S(z)\chi^+\|\}<1,
\end{equation}
which guarantees that the ``interaction'' of the positive and negative ``parts'' of the perturbation $V$  is ``small''. This is to say that it restricts the class of operators considered to those which are not too far from an orthogonal sum of a dissipative and an anti-dissipative ($\im L \leq 0$) operators.

In particular, \cite{Naboko_singular} provides non-restrictive additional conditions such that
$$
N_i(L)\cap N_e (L)=\{0\}, \quad N_i(L)\vee N_e(L) =K
$$
and sharp estimates for the angle between $N_i(L)$ and $N_e (L)$. What's more,
$$
N_-^i(L) \cap N_+^i(L)=\{0\}, \quad N_-^i(L)\vee N_+^i (L)= N_i(L)
$$
with an explicit estimate for the angle between $N_-^i(L)$ and $N_+^i(L)$. Further, $L|_{N_+^i}$ ($L|_{N_-^i}$) is similar to a dissipative (anti-dissipative, respectively) operator with purely singular spectrum.

Dropping the separability condition \eqref{eq:separability} makes the spectral
analysis of $L$ much more involved. The corresponding problems were posed by S.
Naboko in \cite{Naboko_problems}. Most of them are still awaiting resolution,
including the problem of a general spectral resolution of identity for a
non-self-adjoint operator of the class considered here, but some were
successfully tackled in \cite{Naboko_Veselov} by S. Naboko and his student V.
Veselov as well as in subsequent papers of V. Veselov. In particular, the named
paper concerns with an in-depth study of the spectral subspace $N_0^i(L)$,
introduced in \cite{Naboko_problems}.
The main result is formulated {for $V\in \mathfrak{S}_1$} as follows:
$$
\det \Theta(z) = \det \Theta_{L|_{N_+(L)\vee N_-(L)}}(z), \quad N_+(L)\vee N_-(L)=N_-^i(L)\vee N_e(L) \vee  N_+^i(L),
$$
generalising the corresponding result of Gohberg and Kre\u\i n. It shows that the determinant of the characteristic function of $L$ contains no information on the spectral subspace
\begin{equation}\label{eq:Nio}
N_0^i(L):=   K\ominus \{N_-(L^*)\vee N_+(L^*)\}\subset N_i(L),
\end{equation}
i.e., $\det \Theta_{L|_{N_0^i(L)}}(z)\equiv 1.$
Here in notation of Section \ref{sec:smooth} $N_\pm(L^*)=N_\pm^{-iJ}$.

The subspace $N_0^i$ is precisely the ``additional'' spectral subspace corresponding to the real part of the spectrum of $L$ (in particular, it contains the eigenvectors and root vectors corresponding to real values of the spectral parameter), the analytic structure of which has no parallels in the case of dissipative operators. In a nutshell, it appears due to the interaction of the ``incoming'' and ``outgoing'' energy channels in the non-conservative system modelled by $L$.

The r\^ole of $N_0^i$ for the spectral analysis of non-dissipative operators is
further revealed by the following assertion:
$$
N_i(L)\cap N_e(L)\subset N_0^i(L),
$$
i.e., \emph{if} the absolutely continuous and singular subspaces intersect, the
intersection \emph{must} lie in $N_0^i$. It is therefore the presence of $N_0^i$
that ensures that $N_e(L^*)\vee N_i(L^*)\not = K$, which prevents a spectral
decomposition for the
 operator $L^*$.

Sergey had mentioned to us, that he had seven to eight papers worth of further material on the functional model and spectral analysis of non-dissipative operators. Unfortunately, he had never published these results.

\subsection{A functional model based on the Strauss characteristic function}

In contrast to the model theory for contractions associated
with the names of Sz.-Nagy-Foia\c{s} and de~Branges-Rovnyak,
the models of unbounded non-selfadjoint operators are usually
concerned with ``concrete'' operators
arising in applications.
In particular, the functional model
for non-self-adjoint additive
perturbations discussed above
was motivated by the spectral analysis of the Schr\"{o}dinger
operator with a complex potential, see, e.g., \cite{MR0510053, nabokoSpectralComponents,
Naboko_singular, MR1252228}. In fact, Sergey Naboko had reiterated to us on a number of occasions, that his primary concern was the spectral theory of the Schr\"odinger operator, rather than the development of abstract mathematical concepts: the functional model in his view was simply the tool of choice in this area.
More precisely,  the Schr\"odinger
operator~$-\Delta + p(x) + i q(x)$ in ~$L_2(\mathbb R^3)$,
where~$p(x)$, $q(x)$ are real-valued bounded functions
of~$x \in \mathbb R^3$, can be written in the
form~(\ref{eqn:AV}) with the operator~$\alpha$
defined as  $\alpha: f \mapsto
 |q(x)|^{1/2}f$, where $f\in L_2(\mathbb R^3)$.
It is important to note
that all the building blocks of  the
model construction are explicitly given
in terms of the problem at hand.
Indeed, both
the characteristic function~$S(z)$ and the
``spectral maps''~$\mathscr F_\pm$ are
expressed via the non-real part of
the complex potential
(and the operator itself).
The true nature of the problem's
``non-selfadjointess'', i.e., the non-triviality
of the imaginary part of the potential,
is thus faithfully preserved in the
model representation.

The same observation is valid
for other model constructions
of non-selfadjoint operators available in the
literature, see, e.g., \cite{CEKRS} in the present volume for the case of non-self-adjoint extensions of symmetric operators.
%
%
Therefore it becomes increasingly
important to express the
non-selfadjointness of the problem
not in abstract terms (as it is
commonly done in the operator theory),
but rather in terms of the concrete
operator present in
the problem statement.

The standard way to calculate
the characteristic
function of a non-self-adoint operator
is based on the
definition given
by A.~Strauss
in~\cite{Strauss1960}.
For a dissipative operator
it reads as follows
\begin{definition}[\cite{Strauss1960}]
Let ~$L$ be a closed maximal densely defined
dissipative operator on a Hilbert
space~$K$.
The characteristic function of $L$
is a bounded operator-valued analytic
function~$S(z): E \to E_*$,
$z\in\rho(L^*)$, such that
$$
S(z)\Gamma f = \Gamma_* (L^* - zI)^{-1}(L -zI)f,
\quad f \in \dom(L),
$$
where the {\it boundary operators}~$\Gamma$, $\Gamma_*$
are
defined
for $u,v\in \dom(L)$, $u^\prime, v^\prime \in \dom(L^*)$
by the equalities
$$
 (A u, v) - (u, Av )= i(\Gamma u, \Gamma v)_E, \qquad
 ( u^\prime, A^*v^\prime) - (A^*u^\prime, v^\prime )
 = i(\Gamma_* u^\prime, \Gamma_*
v^\prime)_{E_*}
$$
and~$E := \clos\ran(\Gamma)$,
$E_* := \clos\ran(\Gamma_*)$ are Hilbert spaces.
\end{definition}

According to this definition, the concrete form of
the characteristic function of~$L$ depends
on the choice of
boundary operators~$\Gamma$,
$\Gamma_*$.
It is easy to see that for any Hilbert space
isometries~$\pi: E \to E^\prime$,
$\pi_*: E_* \to E_*^\prime$, the
maps~$\pi\Gamma$ and $\pi_*\Gamma_*$
are also boundary operators with
the corresponding characteristic
function~$\pi_* S(z) \pi^* : E^\prime \to E_*^\prime$.
In applications, a suitable definition of the
boundary operators is determined according
to the problem statement itself.
For example, the operator~$\alpha$
of~(\ref{eqn:AVMod}) (the root cause
of the operator's non-selfadointness)
admits the r\^ole
of both~$\Gamma$ and $\Gamma_*$.
Convenient boundary operators 
appear ``naturally'' in the analysis of
non-self-adjoint extensions
of symmetric operators as well.
Once the triple~$\{\Gamma, \Gamma_*, S(z)\}$
is explicitly defined,
the construction of the functional model
  follows the
blueprint of S.~Naboko~\cite{MR573902}.

A further important contribution is contained in
the two recent papers~\cite{BMNW2018,BMNW2022} by B.M. Brown, M. Marletta, S. Naboko, and I. Wood. The authors offer
a model construction carried out
in the abstract setting of {S}trauss'
boundary operators~$\Gamma$,
$\Gamma_*$, resorting to
no specific realisation of them.
This work therefore makes all the steps
of the model construction explicit,
regardless of any particular form of
the characteristic function, the latter to be set based on the
requirements imposed by a concrete application at hand. In particular, this makes it possible to construct a functional model in the case where both the differential expression itself and the boundary conditions are non-self-adjoint, which in our view is especially relevant for topical problems of materials science.


\subsection{Applications of the functional model technique}

Here we list some notable applications of the functional model technique, in which Sergey Naboko was involved, in addition to his work on the spectral analysis of non-self-adjoint Schr\"odinger operators mentioned earlier, see, e.g., \cite{Naboko_singular} and references therein.

1. In  \cite{Romanov_faa,Romanov_indiana} Sergey, together with Yu. Kuperin and R. Romanov, studied the non-self-adjoint single-velocity Boltzmann transport operator. Using the functional model techniques, the absolute continuity of this operator's continuous spectrum was proved; the similarity problem of the absolutely continuous ``part'' of the operator to a self-adjoint one was fully settled, and the existence of a spectral singularity at zero ascertained for a singular set of multiplication coefficients.

2. In \cite{Romanov1,Romanov2}, together with R. Romanov, Sergey Naboko analysed the impact of spectral singularities on the asymtotic behaviour of the group of exponentials, generated by a maximal dissipative operator $L$. It was shown that this asymptotics allows one to recover the orders and locations of spectral singularities in the case, where their number is finite and they are of a finite power order.

3. In \cite{Kiselev_faa,Kiselev_arkiv}, for a non-dissipative trace class perturbation $L$ of a self-adjoint operator on $K$ such that $N_0^i(L)$ coincides with the Hilbert space $K$, a generalisation of the Caley identity was obtained in the following form: there exists an outer in the upper half-plane $\mathbb
C_+$ uniformly bounded scalar analytic function $\gamma(\lambda)$ such that
$
w-\lim_{\varepsilon\downarrow 0} \gamma(L+i\varepsilon) =0.
$
A generalisation of this result was further obtained to the case of relative trace class perturbations.

4. In \cite{Kiselev_jcam}, the so-called matrix model was introduced and studied in some detail, i.e., a rank two non-dissipative additive perturbation $L$ in $K$ of a self-adjoint operator under the assumption that $K=N_0^i(L)$. This model represents the simplest possible case of a non-dissipative operator which exhibits the properties not found in any dissipative one; despite its seeming simplicity, it already includes the main analytic obstacles found in the general case. It has to be noted that this model was the favourite sandbox of Sergey; unfortunately, many results obtained by him, up to and including a von Neumann type estimate in BMO classes for functions of the operator $L$, have never been published.



\end{document}